%% file: power-bound-generalization.tex
\documentclass[12pt, hidelinks]{article}
\usepackage[utf8]{inputenc}
\usepackage{amsmath, amssymb, url, float, graphicx}
\usepackage{fullpage}
\usepackage[small,bf]{caption}

\DeclareMathOperator{\minimize}{minimize}

\renewcommand\phi\varphi

\input defs.tex

\title{A Note on Generalizing Power Bounds for Physical Design}
\author{Guillermo Angeris}
\date{August 2022}

\begin{document} 
\maketitle 

\begin{abstract}
    In this note we show how to construct a number of very general nonconvex
    quadratic inequalities for a variety of physics equations appearing in
    physical design problems. These nonconvex quadratic inequalities can then
    be used to construct bounds on physical design problems where the objective
    is a quadratic or a ratio of quadratics. We show that the quadratic
    inequalities and the original physics equations are equivalent under a
    technical condition that holds in many practical cases which is easy to
    computationally (and, in some cases, manually) verify.
\end{abstract}

\section*{Introduction}
The physical design problem appears in many contexts including in photonic
design, antenna design, horn design, among many others. In general, almost all
formulations have very similar structure, and, recently, there has been a large
amount of work on fast algorithms for approximately solving the physical design
problem~\cite{luNanophotonicComputationalDesign2013,menRobustTopologyOptimization2014,
    piggottInverseDesignDemonstration2015,moleskyInverseDesignNanophotonics2018,
doryInversedesignedDiamondPhotonics2019, kesiciMethodAntennaImpedance2020,
angerisConvexRestrictionsPhysical2021}, which is generally NP-hard to solve,
and bounding the optimal value of such
problems~\cite{millerFundamentalLimitsOptical2016,
    angerisComputationalBoundsPhotonic2019,
    moleskyTOperatorLimitsElectromagnetic2020,
    gustafssonUpperBoundsAbsorption2020,
    trivediBoundsScatteringAbsorptionless2020,
    kuangComputationalBoundsLight2020,
    moleskyGlobalOperatorBounds2020,
    zhaoMinimumDielectricResonatorMode2020},
which in general cannot be efficiently found (except, perhaps, in
some special cases, \cf~\cite{moleskySionMinimaxTheorem2021}). In this
note, we show how to reduce a large family of these problems to a problem
over a number of nonconvex quadratic constraints. This family includes a
number of problems which exhibit `nonlocal scattering'; \ie, problems where
some of the design parameters can affect many, if not all, points in the
domain. These nonconvex quadratic constraints can then be easily relaxed to
obtain efficiently-computable lower bounds for objectives which are also
quadratics or ratios of quadratics. As this is a short note, we assume a
reasonable amount of familiarity with physical design problems. For some
overviews of this problem and further references, see,
\eg,~\cite{angerisHeuristicMethodsPerformance2021,
chaoPhysicalLimitsElectromagnetism2022}.

\section{The physical design problem}
We will define the physical design problem in this section and then show an equivalent
optimization problem over only the fields $z$. We will then see, in the next section,
how to find the dual to this new problem to get lower bounds on the best possible solution.

\paragraph{Physics equation.} 
In general, the designer begins with some \emph{physics equation}, which we
will write as
\begin{equation}\label{eq:physics}
    A(\theta) z= b.
\end{equation}
Here, $A(\theta) \in \reals^{m\times n}$ is the \emph{physics matrix},
$b \in \reals^m$
is the \emph{excitation}, while $z \in \reals^n$ is the \emph{field}, all
of which are parametrized by some set of \emph{design parameters},
$\theta \in \reals^d$. We will assume, as is the case in many instances of
physical design, that the parameters enter in a very specific way:
\[
A(\theta) = A_0 + \sum_{i=1}^d \theta_i A_i.
\]
Here, the matrices satisfy $A_i \in \reals^{m \times n}$. We call this specific
parametrization \emph{affine}, as $A$ is an affine function of the design
parameters, $\theta$. Note that, unlike many of the currently-known bounds for
physical design problems~\cite{angerisComputationalBoundsPhotonic2019,
    moleskyTOperatorLimitsElectromagnetic2020,
    gustafssonUpperBoundsAbsorption2020,
    trivediBoundsScatteringAbsorptionless2020,
    kuangComputationalBoundsLight2020,
    moleskyGlobalOperatorBounds2020,
    zhaoMinimumDielectricResonatorMode2020},
we do not assume that the matrices
$A_i$ are outer products of the unit basis vectors ($A_i = e_ie_i^T$) or
even diagonal.

\paragraph{Parameter constraints.} We will assume that the parameters $\theta$
are constrained to lie in some interval.
From~\cite[\S2.2]{angerisHeuristicMethodsPerformance2021} we know that, in the
special case of the affine parameterization, we may generally assume that
$\theta$ lies in the following interval
\[
-\ones \le \theta \le \ones.
\]
(Or, in the case of Boolean constraints, that $\theta$ lies in the unit
hypercube, $\theta \in \{\pm 1\}^d$; we will see extensions to this case
later.) 

\paragraph{Objective function and problem.} We will assume that the designer
wishes to minimize some objective function $f: \reals^n \to \reals$ that need
not be convex, which depends only on the field $z$. Putting this all together,
we have the following (nonconvex) optimization problem, which we will call the
\emph{physical design problem}:
\begin{equation}\label{eq:main}
\begin{aligned}
	& \text{minimize} && f(z)\\
	& \text{subject to} && A(\theta)z = b\\
	&&& -\ones \le \theta \le \ones,
\end{aligned}
\end{equation}
with variables $z \in \reals^n$ and $\theta \in \reals^d$, where $A$ is
affine in the parameters.

\subsection{Equivalent constraints}\label{sec:equivalence}
To start, we will show a simple equivalent characterization of when two real
vectors $x$ and $y$ are collinear, $x = \alpha y$, with scale factor $|\alpha|
\le 1$. We will then show how this characterization can be used in rewriting the
optimization problem~\eqref{eq:main}, which depends on both the fields $z$ and
the parameters $\theta$, to an equivalent problem that depends only on the
fields $z$; \ie, we will show how to eliminate the parameters $\theta$ from
this problem and give a set of conditions for when this procedure is
tight.

\paragraph{Equivalent characterization.} Given two vectors $x, y \in \reals^n$,
we will show that $x = \alpha y$ with $-1 \le \alpha \le 1$ if, and only if
\begin{equation}\label{eq:psd-condition}
    x^TNx \le y^TNy, ~ \text{for all} ~ N \in \symm_+^n,
\end{equation}
where $N \in \symm_+^n$ is the set of positive semidefinite matrices. The
forward implication is fairly simple: note that if $x = \alpha y$ then, for any
$N \in \symm^n_+$,
\[
x^TNx = \alpha^2 y^TNy \le y^TNy,
\]
since $\alpha^2 \le 1$ and $y^TNy \ge 0$ as $N$ is positive semidefinite.

The backward implication, if $x^TNx \le y^TNy$ for all PSD matrices $N$, then
$x = \alpha y$ is slightly trickier. To see this, since $x^TNx \le y^TNy$ is true for any PSD matrix $N$ by assumption,
then we will choose $N = vv^T$ where
\[
v = (y^Ty) x - (x^Ty)y.
\]
This yields
\[
(v^Tx)^2 = x^TNx \le y^TNy = (v^Ty)^2 = ((y^Ty) x^Ty - (x^Ty)y^Ty)^2 = 0.
\]
The left hand side of this inequality then satisfies
\[
(v^Tx)^2 = ((y^Ty)(x^Tx) - (x^Ty)^2)^2 \le 0,
\]
where the inequality follows from the previous. In other words,
\[
(y^Ty)(x^Tx) - (x^Ty)^2 = 0,
\]
or that
\[
(x^Ty)^2 = \|x\|_2^2\|y\|_2^2.
\]
By Cauchy--Schwarz, this happens if, and only if, $x = \alpha y$ for some
$\alpha \in \reals$. Finally, because $x^TNx \le y^TNy$ for any PSD $N$, we
will choose $N = I$ to get
\[
\alpha^2 y^Ty = x^Tx \le y^Ty,
\]
so $\alpha^2 \le 1$, which happens if, and only if, $-1 \le \alpha \le 1$, as
required. Putting these two together, we get the final result that $x^TNx \le
y^TNy$ for all $N \in \symm^n_+$ if, and only if, $x = \alpha y$ for some $-1
\le \alpha \le 1$. (In fact, note that it suffices for $N$ to be in the set of
rank-one symmetric PSD matrices and the identity. This suggests a possibly
memory-efficient way of dealing with these inequalities in practice.)

\paragraph{Constructing inequalities.} To construct some inequalities, we will now
introduce the family of matrices $P_i \in \reals^{m_i \times m}$ for
$i=0, \dots, d$. We assume that the matrices $P_i$ for $i=1, \dots, d$ have the following property:
\begin{equation}\label{eq:p-conditions}
    P_iA_j = 0,~ \text{whenever} ~ i\ne j, ~ \text{for} ~ i, j= 1, \dots, d,
\end{equation}
while $P_0 \in \reals^{m_0 \times m}$ satisfies:
\begin{equation}\label{eq:zero-p}
    P_0 A_i = 0, ~ \text{for} ~ i=1, \dots, d.
\end{equation}
(In some cases we might have $m_0 = 0$, so, for convenience, we will say that
this equality is then trivially satisfied.) These matrices then have the
property that, for any $\theta$ and $z$ satisfying the physics
equation~\eqref{eq:physics}, when $i=1, \dots, d$,
\[
P_iA(\theta)z = P_i(A_0 + \theta_iA_i)z = P_ib.
\]
In other words, the matrix $P_i$ `picks out' the $i$th design parameter.
Rearranging slightly, we then have that
\[
P_i(b - A_0z) = \theta_i P_iA_iz,
\]
must be satisfied for $-1 \le \theta_i \le 1$ for $i=1, \dots, d$, and
\[
P_0A(\theta)z = P_0A_0z = P_0b.
\]
From condition~\eqref{eq:psd-condition}, we know the first happens if, and only if,
\begin{equation}\label{eq:psd-inequality}
    (b - A_0z)^TP_i^TN_iP_i(b - A_0z) \le z^TA_i^TP_i^TN_iP_iA_iz, ~ \text{for all} ~ N_i \in \symm_+^{m_i},
\end{equation}
for each $i=1, \dots, d$, which is a family of (potentially nonconvex)
quadratic inequalities over $z$, while the second is an affine constraint on
$z$. The remaining question is: when is this set of inequalities tight? In
other words, if $z$ satisfies these inequalities, when do design parameters
$\theta \in \reals^d$ with $-\ones \le \theta \le \ones$ exist such that $z$
and $\theta$ satisfy the physics equation~\eqref{eq:physics}?

\paragraph{Tightness.} As mentioned before, note that~\eqref{eq:psd-inequality}
does not always lead to a tight set of inequalities. For example, setting all
of the $P_i = 0$ for $i=1, \dots, d$, with, say $m_i = 1$,
satisfies~\eqref{eq:p-conditions}, but leads to a trivial set of inequality
constraints that are satisfied for any $z$. One simple set of conditions is the
existence of matrices $M_i \in \reals^{m \times m_i}$ such that
\begin{equation}\label{eq:m-condition}
    \sum_{i=0}^d M_iP_i = I.
\end{equation}
This works because any $z$ that satisfies~\eqref{eq:psd-inequality} also satisfies
\[
P_iA(\theta)z = P_ib,
\]
from the previous discussion. Multiplying both sides of this equation by $M_i$ on
the left hand side and summing over $i$ gives
\[
\sum_{i=0}^d M_iP_iA(\theta)z = A(\theta)z = \sum_{i=0}^d M_iP_ib = b,
\]
which is just the original physics equation in problem~\eqref{eq:main}.

\subsection{A (computational) sufficient condition}\label{sec:sufficient-condition}
While the above discussion is potentially useful as a theoretical tool, it is
often not clear how to generate the matrices $P_i$ (and, in turn, to know when
the matrices $M_i$ exist for this choice of $P_i$). Below we give a simple
sufficient condition which, given matrices $A_i$ satisfying this condition, can
be used to construct $P_i$ and $M_i$ which satisfy
conditions~\eqref{eq:p-conditions}, \eqref{eq:zero-p},
and~\eqref{eq:m-condition} for $i=0, \dots, d$.

\paragraph{Tightness result.} There exists a fairly general check for tightness
that can be easily (computationally) verified. In particular, we will assume
the matrices $A_i$ are written in the form:
\[
A_i = U_iV_i^T,
\]
where $U_i \in \reals^{m \times m_i}$ and $V_i \in \reals^{n \times m_i}$ for
$i=1, \dots, d$. (This can be done by computing, say, the reduced singular
value decomposition for each matrix $A_i$. Often, though, this decomposition is
known in practice, as we will see in the examples.) Now, define the matrix
\[
U = \begin{bmatrix}
    U_1 & U_2 & \dots & U_d
\end{bmatrix}.
\]
The final condition is that, if $U$ is full column rank (\ie, the columns of
$U$ are linearly independent) then we can find matrices $P_i \in \reals^{m_i
\times m}$ satisfying the required conditions, and corresponding matrices $M_i
\in \reals^{m \times m_i}$ satisfying~\eqref{eq:m-condition}. Note that we have
no additional conditions on the matrices $V_i$ (in the definition of $A_i$) nor
any conditions on $A_0$ or $b$, apart from their dimensions.

\paragraph{Proof.} To see this, let $U_0$ be the basis completion of $U$ such
that $\tilde U = \begin{bmatrix}U_0 & U\end{bmatrix}$ is square and invertible.
(From basic linear algebra, we know a $U_0$ exists for any `tall' matrix and
that $\tilde U$ is invertible since it is full column rank.) Now we can set
\begin{equation}\label{eq:p-def}
    \begin{bmatrix}
        P_0\\
        P_1\\
        P_2\\
        \vdots\\
        P_d
    \end{bmatrix} = \tilde U^{-1},
\end{equation}
where $P_i \in \reals^{m_i \times m}$ for $i=0, \dots, d$, where the $m_i$ come
from the size of the $U_i$ defined above, while $m_0 = m - \sum_{i=1}^d m_i \ge
0$. (The inequality here comes from the assumption that $U$ has full column
rank.) With this, we have, using the fact that $\tilde U^{-1} \tilde U = I$,
\[
P_iA_j = (P_iU_j)V_j^T = 0,
\]
if $i\ne j$ for every $i, j=1, \dots, d$, while $P_0 A_j = 0$. (We also have,
of course, that $P_i U_i = I$.) The tightness condition is then met by setting
$M_i = U_i$ for $i=0, \dots, d$ since
\[
    \sum_{i=0}^d M_iP_i = \sum_{i=0}^d U_iP_i = \tilde U \tilde U^{-1} = I,
\]
as required.

\subsection{Examples}
There are many important special cases for which the condition above can be
both easily verified by hand and the matrices $P_i$ (and corresponding $M_i$)
are also easily intuited. In these examples, the procedure above yields the
same results as computing the matrices by hand.

\paragraph{Multi-scenario design.} One classic example in the case of physical
design is the multi-scenario case where the matrices $A_i$ are diagonal with
nonoverlapping nonzero entries, \ie, $A_iA_j = 0$ for $i, j=1, \dots, d$ and $i
\ne j$. If $m_i$ is the number of nonzero entries of $A_i$ then we can write
$A_i = U_iV_i^T$ for $U_i \in \{0, 1\}^{m \times m_i}$ and $V_i \in \reals^{m
\times m_i}$ for $i=1, \dots, d$. We can then set
\[
    \tilde P_0 = I - \sum_{i=1}^d U_iU_i^T, \quad P_i = U_i^T, ~~i=1, \dots, d.
\]
and let $P_0\in \reals^{m_0 \times m}$ be $\tilde P_0$ with zero-rows removed
and $m_0 = m - \sum_{i=1}^d m_i$. (We can also view $P_0$ as the matrix such
that $P_0^T$ has all of the unit vectors not appearing in the columns of the
matrices $U_1, \dots, U_d$.) Setting $M_i = P_i^T$ for $i=0, \dots, d$ then
satisfies all required conditions. We note that this common special case was
first shown in~\cite{shimFundamentalLimitsMultifunctional2021} and the
construction presented here results in an identical formulation for the dual
(presented in a later section).

\paragraph{Rank-one matrices.} Another important special case is when the
matrices $A_i = u_iv_i^T$ are rank-one matrices for $i=1, \dots, d$. In this
case, the conditions above state that $U$ is a matrix whose columns are the
vectors $u_i$ and that the vectors $u_i$ must be linearly independent. Setting
the $P_i$ matrices as in~\eqref{eq:p-def} leads to $d$ row vectors $P_i = p_i^T
\in \reals^{1\times m}$ which satisfy
\[
    p_i^Tu_j = \begin{cases}
        1 & i = j\\
        0 & \text{otherwise},
    \end{cases}
\]
for $i, j=1, \dots, d$ with $i \ne j$, while $P_0$ is a basis
for the subspace orthogonal to the space generated by the $\{u_i\}$.
Assuming that $P_0$ is normalized to $P_0^TP_0 = I$ then setting $M_0 = P_0^T$
and $M_i = u_i$ for $i=1, \dots, d$ gives the desired conditions.

\paragraph{Discussion.} These examples let us interpret
equations~\eqref{eq:p-conditions} and~\eqref{eq:m-condition} as a condition
that the matrices $A_i$ are not `too correlated' in the sense that their left
singular vectors don't overlap `too much'.

\paragraph{Final problem.} Replacing the physics equation and the
constraints over $\theta$ with the set of
inequalities~\eqref{eq:psd-inequality} then gives the final problem:
\begin{equation}\label{eq:final-problem}
\begin{aligned}
	& \text{minimize} && f(z)\\
	& \text{subject to} && z \in S_i, \quad i=1, \dots, d\\
    &&& P_0A_0z = P_0b,
\end{aligned}
\end{equation}
with variable $z \in \reals^n$, where $S_i \subseteq \reals^n$ is defined as
\[
    S_i = \{z \in \reals^n \mid (b - A_0z)^TP_i^TN_iP_i(b - A_0z) \le z^TA_i^TP_i^TN_iP_iA_iz, ~ \text{for all} ~ N_i \in \symm^{m_i}_+\},
\]
for $i=1, \dots, d$. In other words, $S_i$ is the set of fields $z$ consistent
with condition~\eqref{eq:psd-inequality}. Note that, as written, $z \in S_i$ denotes an infinite
family of constraints (over all possible positive semidefinite matrices of
dimension $m_i$) which need not even be convex. For future reference, we will
denote the optimal value of this problem $p^\star$.

\paragraph{Extensions.} There are some simple extensions to this formulation that
are similarly useful. For example, if we constrain $\theta \in \{\pm1\}^d$, \ie,
require the parameters $\theta$ to be Boolean, the resulting problem is
\[
\begin{aligned}
	& \text{minimize} && f(z)\\
	& \text{subject to} && z \in \tilde S_i, \quad i= 1, \dots, d\\
    &&& P_0A_0z = P_0b,
\end{aligned}
\]
where
\[
    \tilde S_i = \{ z \in \reals^n \mid (b - A_0z)^TP_i^TN_iP_i(b - A_0z) \le z^TA_i^TP_i^TN_iP_iA_iz, ~ \text{for all} ~ N_i \in \symm^{m_i}\}.
\]
Note that $N_i$ is now allowed to range over all symmetric matrices, not just
the positive semidefinite ones.
To see this equivalence, note that we know
\[
x^TNx \le y^TNy, ~ \text{for all} ~ N \in \symm^n,
\]
implies that $x = \alpha y$ for $\alpha^2 \le 1$ from the same argument
as the characterization in~\S\ref{sec:equivalence}. If we also set $N = -I$, we have that
\[
x^Tx \ge y^Ty,
\]
so $\alpha^2 \ge 1$ and therefore $\alpha^2 = 1$, so $\alpha \in\{\pm 1\}$. The
converse (that $x = \alpha y$ implies the inequality above when $\alpha \in
\{\pm 1\}$) is easily verified. Performing the same steps as before yields this
problem. We may also write the sets $\tilde S_i$ as
\[
    \tilde S_i = \{z \in \reals^n \mid (b - A_0z)^TP_i^TN_iP_i(b - A_0z) = z^TA_i^TP_i^TN_iP_iA_iz, ~ \text{for all} ~ N_i \in \symm^{m_i}_+ \}
\]
where we have replaced the inequality with an equality, but only allow $N_i$ to
range across the positive semidefinite matrices.

\section{Dual problem}\label{sec:dual-problem}
We will now show how to compute a dual problem of~\eqref{eq:final-problem}. 

\paragraph{Rewriting.} One simple way of rewriting the problem~\eqref{eq:final-problem}
is to include the constraints over the sets $S_i$ as an indicator function in the objective
$I_i: \reals^n \to \reals \cup \{+\infty\}$, defined
\[
I_i(z) = \begin{cases}
    0 & z \in S_i\\
    +\infty & \text{otherwise},
\end{cases}
\]
for $i=1, \dots, d$ and, for convenience, we will define
\[
    I_0(z) = \sup_\nu \nu^T(P_0A_0z - P_0b).
\]
which we note is zero if $P_0A_0z - P_0b = 0$ and is positive infinity
otherwise. (In other words, $I_0$ is an indicator function for the affine
constraint.) Then, the following problem is equivalent to~\eqref{eq:final-problem}:
\[
\begin{aligned}
    & \text{minimize} && f(z) + \sum_{i=0}^d I_i(z)\\
\end{aligned}
\]
with variable $z$.

\paragraph{Indicator function.} We will now see that the indicators $I_i$ for
$i=1, \dots, d$ have a convenient form, similar to that of $I_0$. In
particular, we will show that we can write
\[
I_i(z) = \sup_{N_i \ge 0} \left((b - A_0z)^TP_i^TN_iP_i(b - A_0z) - z^TA_i^TP_i^TN_iP_iA_iz\right),
\]
for each $i=1, \dots, d$. To see this, consider first that if $z \in S_i$, then
\[
    (b - A_0z)^TP_i^TN_iP_i(b - A_0z) \le z^TA_i^TP_i^TN_iP_iA_iz
\]
for every $N_i \in \symm^{m_i}_+$, by definition of $S_i$. So, using this
definition of $I_i$ we would have that
\[
    I_i(z) \le 0.
\]
Picking $N_i = 0$ in the supremum above saturates the inequality, so $I_i(z) =
0$ if $z \in S_i$.

On the other hand, if $z \not \in S_i$ then there exists some $i$ and
$N_i \ge 0$ with
\[
(b - A_0z)^TP_i^TN_iP_i(b - A_0z) > z^TA_i^TP_i^TN_iP_iA_iz.
\]
Letting $t > 0$ be any positive value, then $N_i$ is PSD implies that $tN_i$ is also PSD,
so we have that
\[
I_i(z) \ge t((b - A_0z)^TP_i^TN_iP_i(b - A_0z) - z^TA_i^TP_i^TN_iP_iA_iz).
\]
Sending $t \to \infty$ implies that $I_i(z) = \infty$ when $z \not \in S_i$ as required.

\paragraph{Rewritten problem.} We can then rewrite problem~\eqref{eq:final-problem}
as a saddle point problem:
\begin{multline*}
	\minimize \sup_{N \ge 0, \nu}\bigg(f(z) + \sum_{i=1}^d \left((b - A_0z)^TP_i^TN_iP_i(b - A_0z) - z^TA_i^TP_i^TN_iP_iA_iz\right) \\+ \nu^T(P_0A_0z - P_0b)\bigg),
\end{multline*}
with variable $z \in \reals^n$. We have written $N \ge 0$ as a shorthand for
$N_i \in \symm_+^{m_i}$ for $i=1, \dots, d$, and have pulled the supremum
outside of the sum. We can view the function
\[
L(z, N, \nu) = f(z) + \sum_{i=1}^d \left((b - A_0z)^TP_i^TN_iP_i(b - A_0z) - z^TA_i^TP_i^TN_iP_iA_iz\right) + \nu^T(P_0A_0z - P_0b),
\]
as a Lagrangian for problem~\eqref{eq:final-problem} in that a solution
to the original problem is a saddle point of this function, maximizing over $N$
and then minimizing over $z$. In particular, from the previous discussion,
we know that
\[
p^\star = \inf_z \sup_{N \ge 0, \nu} L(z, N, \nu),
\]
where $p^\star$ is the optimal value of problem~\eqref{eq:final-problem}.
In fact, this suggests a reasonable heuristic for~\eqref{eq:final-problem}
is to use a saddle-point solver. We suspect this might be relatively
efficient as computing $L$ requires only vector-matrix multiplies
and computing the gradient of $L$ with respect to $z$ or $N_i$ for $i=1, \dots, d$
requires only matrix-matrix multiplications of relatively small size.

\paragraph{Bounds.} We can get a simple bound on the optimal value of
problem~\eqref{eq:final-problem}, $p^\star$, by swapping the infimum and the supremum.
To see this, we define the dual function:
\[
g(N, \nu) = \inf_z L(z, N, \nu) \le L(z, N, \nu),
\]
for all $N \ge 0$ and $\nu \in \reals^{m_0}$, so
\[
\sup_{N \ge 0, \nu} g(N, \nu) \le \sup_{N \ge 0, \nu} L(z, N, \nu).
\]
Taking the infimum over $z$ of the right hand side gives
\[
\sup_{N \ge 0, \nu} g(N, \nu) \le \inf_z\sup_{N \ge 0, \nu} L(z, N, \nu) = p^\star.
\]
We will denote this lower bound as
\[
d^\star = \sup_{N \ge 0, \nu} g(N, \nu).
\]
From before, we know that $d^\star$ need not be easy to evaluate. On the other
hand, since $g$ is defined as the infimum of a family of affine functions of 
the matrices $N_i$, it is always a concave function. Whenever $g(N, \nu)$ is easy
to evaluate, it is almost always the case that we can efficiently solve for
$d^\star$ as finding the value is just maximizing a concave function over the 
PSD matrices $N_1, \dots, N_d$ and the vector $\nu \in \reals^{m_0}$.

\paragraph{Evaluating the dual function for a quadratic objective.}
In general, $g$ is not easy to evaluate. On the other hand, in some special
cases, we can give closed-form solutions for $g(N, \nu)$. For example, in the special
case that the objective function $f$ for problem~\eqref{eq:final-problem} is a quadratic,
\ie,
\[
f(z) = \frac12 z^TQz + q^Tz + r,
\]
then $L$ is a quadratic function of $z$ (holding all other variables constant)
which is always easy to minimize. In particular, we can write
\[
L(z, N, \nu) = \frac{1}{2}z^TT(N)z + u(N, \nu)^Tz + v(N, \nu),
\]
where
\[
T(N) = Q + 2A_0^T\left(\sum_{i=1}^d P_i^TN_iP_i\right)A_0 - 2\sum_{i=1}^d A_i^TP_i^TN_iP_iA_i,
\]
and
\[
u(N, \nu) = q - 2A_0^T\left(\sum_{i=1}^d P_i^TN_iP_i\right)b + A_0^TP_0^T\nu, \qquad v(N, \nu) = r + b^T\left(\sum_{i=1}^d P_i^TN_iP_i\right)b - \nu^T P_0b.
\]
Using this, we then have that
\[
g(N, \nu) = \inf_z L(z, N, \nu)= v(N, \nu) - \frac12 u(N, \nu)^TT(N)^+u(N, \nu),
\]
whenever $T(N) \ge 0$ and $u(N, \nu) \perp \nullspace(T(N))$, \ie, $u(N, \nu)$ is
orthogonal to the nullspace of $T(N)$, and where $T(N)^+$ is the Moore--Penrose
pseudoinverse of $T(N)$. Otherwise, $g(N, \nu) = -\infty$.

\paragraph{Dual problem.}
Maximizing the dual function $g$ can be written as the (almost) standard form
semidefinite problem, in inequality form:
\[
\begin{aligned}
    & \text{maximize} && v(N, \nu) - \frac12 t\\
    & \text{subject to} && \begin{bmatrix}
        T(N) & u(N, \nu)\\
        u(N, \nu)^T & t
    \end{bmatrix} \ge 0,
\end{aligned}
\]
with variables $N = (N_1, \dots, N_d)$ with $N_i \in \symm^{m_i}_+$ for $i = 1,
\dots, d$, $\nu \in \reals^{m_0}$ and $t \in \reals$. This transformation is
standard and a simple proof of equivalence follows from the application of a
Schur complement~\cite[\S A.5.5]{cvxbook}.

\paragraph{Extensions.} In the case that $f$ is a ratio of quadratics or has
additional quadratic constraints added as indicator functions, we can apply a
similar transformation as~\cite{angerisBoundsEfficiencyMetrics2022}
which leads to an additional quadratic constraint that is easily included in
the dual formulation.

\section*{Acknowledgements}
The author would like to thank Sean Molesky for comments and suggestions and
Theo Diamandis for, quoting from written communication, ``moral support,''
along with a number of helpful suggestions, edits, and conversations.

\section{Addendum (March 2026)}
While doing other work, I was, for some reason, reminded of this note and
decided to give GPT-5.4 and also Pro, separately, the statement of the key
``lemma''. (That is, if the result can even be called that.) The idea was to
test it and see if it would one-shot the solution, which I assumed would just
recover the original result.

To my surprise, both Thinking (and, in turn, Pro) came back with a slightly
stronger statement and proof that removes the sufficient condition above. My
write-up of GPT-5.4's provided proof with additional commentary is below. 

\subsection{An improved characterization}
As before, we begin with the bilinear system of equations~\eqref{eq:physics}:
\[
    A(\theta)z = b.
\]
As previously done in the paper, we will attempt to eliminate the design
parameters $\theta$ from this equation to get a (potentially infinite) set of
quadratic inequalities over the fields $z$ alone. I specifically asked
GPT-5.4 to give a family of quadratic inequalities that are equivalent to
this original system of equations: that is, if $z$ satisfies the family of
quadratics, then there exists $\theta \in [-1, 1]^d$ such that $A(\theta)z = b$,
and vice versa.

\paragraph{Scalarization.}  Note that this system of equations is true if, and
only if,
\[
    y^TA(\theta)z = y^Tb \quad \text{for all $y \in \reals^m$}.
\]
Equivalently, using the fact that, by assumption, the physics matrix $A(\theta)$
is affine in the design parameters $\theta \in \reals^d$, that is
\[
    A(\theta) = A_0 + \sum_{i=1}^d \theta_i A_i,
\]
means that the original equation~\eqref{eq:physics} is true if, and only if,
\begin{equation}\label{eq:scalarized}
    y^T(b - A_0z) = \sum_{i=1}^d \theta_i (y^TA_iz) \quad \text{for all $y \in \reals^m$}.
\end{equation}

\paragraph{Eliminating the design parameters.} One simple observation is to note
that, for any $x \in \reals^d$ and $\alpha \in \reals$, we have that there exists
$\theta \in [-1, 1]^d$ such that
\[
    \alpha = \sum_{i=1}^d \theta_i x_i,
\]
if, and only if,
\[
    |\alpha| \le \sum_{i=1}^d |x_i|.
\]
Applying this to the right hand side of~\eqref{eq:scalarized} implies that, for
some field $z$, there exists $\theta \in [-1, 1]^d$ such
that~\eqref{eq:scalarized} is true if, and only if,
\begin{equation}\label{eq:abs-inequality}
    |y^T(b - A_0z)| \le \sum_{i=1}^d |y^TA_iz| \quad \text{for all $y \in \reals^m$}.
\end{equation}
In other words, there exists $\theta \in [-1, 1]^d$ such that $A(\theta)z = b$ if, and only if
the above inequality is satisfied for all $y \in \reals^m$.

\paragraph{Quadratic inequalities.} As stated, the above inequality is not easy
to work with. On the other hand, we can turn this family of inequalities over
absolute values into a family of quadratic inequalities by noting that there is
a simple `dual' characterization of the squared-sum of absolute values. For any
$x \in \reals^d$, we have that
\[
\left(\sum_{i=1}^d |x_i|\right)^2 \le \sum_{i=1}^d \frac{x_i^2}{w_i}, \quad \text{for all $w \ge 0$ with $\ones^T w = 1$},
\]
with equality exactly when $w_j = |x_j| / \sum_{i=1}^d |x_i|$, where we define
$0/0 = 0$ and $\alpha/0 = \infty$ for any $\alpha > 0$. (To see this, apply
Cauchy--Schwarz on the vectors $(|x_i|/\sqrt{w_i})_i$ and $(\sqrt{w_i})_i$.)
Introducing a new variable, $w_i$ for each term in the sum on the right hand
side of~\eqref{eq:abs-inequality}, we can rewrite the inequality in the
following way, by squaring both sides:
\begin{equation}\label{eq:new-inequality}
    (y^T(b - A_0z))^2 \le \sum_{i=1}^d \frac{(y^TA_iz)^2}{w_i} \quad \text{for all $y \in \reals^m$ and $w \ge 0$ with $\ones^T w = 1$}.
\end{equation}
It is easy to see this is a (nonconvex) quadratic inequality over $z$, but,
unlike the original quadratic inequalities, the equivalence has no conditions on
the matrices $A_0, \dots, A_d$.

This lets us write the final (equivalent) problem as the following:
\begin{equation}\label{eq:final-problem-2}
\begin{aligned}
	& \text{minimize} && f(z)\\
	& \text{subject to} && (y^T(b - A_0z))^2 \le \sum_{i=1}^d \frac{(y^TA_iz)^2}{w_i} ~ \text{for all $y \in \reals^m$ and $w \ge 0$ with $\ones^T w = 1$},
\end{aligned}
\end{equation}
over the variable $z \in \reals^n$.

\paragraph{Discussion.} Note that the above inequalities are at least as strong
as the previous conditions.  In particular, they are equivalent when the $A_i$
satisfy the sufficient conditions given in~\S\ref{sec:sufficient-condition}, but
they also apply in more general cases. To see this, let $P_i$ for $i=0, \dots,
d$ be matrices satisfying~\eqref{eq:p-conditions} and~\eqref{eq:zero-p}, then
note that we can recover the previous inequalities

\subsection{Getting a lower bound}
Unfortunately, after this, I was not able to reasonably prompt GPT-5.4 to give
me a simple program (as in the original note) nor a nice reduction to the
original set of inequalities under the sufficient conditions above. Below I show
a simple way of building a lower bound that resembles the original dual problem
and inequalities, and show how to recover the original dual problem under the
sufficient conditions of~\S\ref{sec:sufficient-condition}.

\paragraph{Rewriting the family.} One simple observation is that we can take nonnegative
(conic) combinations of the new family of inequalities~\eqref{eq:new-inequality}.
In particular, let $\lambda(y) \ge 0$ be a family of nonnegative weights
with finite support (\ie, $\lambda(y) > 0$ for only finitely many $y$), then
any $z$ that satisfies~\eqref{eq:new-inequality} also satisfies
\begin{equation}\label{eq:conic-combination}
\sum_y \lambda(y) (y^T(b - A_0z))^2 \le \sum_y \lambda(y) \sum_{i=1}^d \frac{(y^TA_iz)^2}{w_i}, ~~ \text{for all $w \ge 0$ with $\ones^T w = 1$}.
\end{equation}
Note that since any positive semidefinite matrix $N \in \symm_+^m$ can be
decomposed as $N = \sum_{i=1}^m \sigma_i v_i v_i^T$ for some $\sigma_i \ge 0$
and $v_i \in \reals^m$ where $i=1, \dots, m$, then any $z$ that
satisfies~\eqref{eq:conic-combination} over all $\lambda(y) \ge 0$ with finite
support also satisfies
\[
(b - A_0z)^T N (b - A_0z) \le \sum_{i=1}^d \frac{z^TA_i^T N A_iz}{w_i}, ~~ \text{for all $N \in \symm_+^m$ and $w \ge 0$ with $\ones^Tw = 1$}.
\]
This follows from~\eqref{eq:conic-combination} for any $N = \sum_{i=1}^{m}
\sigma_i v_i v_i^T$ by setting $\lambda(y) = \sigma_i$ whenever $y = v_i$ for
some $i$ and $\lambda(y) = 0$ otherwise.  Indeed, since $yy^T \in \symm_+^m$ is
positive semidefinite, then both families of inequalities are equivalent in that
any $z$ that satisfies the first must also satisfy the second and vice versa.

\paragraph{On tightness.} In general, the family of
inequalities~\eqref{eq:conic-combination} is true for any measure, so we may
take $\lambda$ to be a measure over $\reals^m$ and we can replace the sum with
an integral. On the other hand, the finite support case suffices since the
integrand is monotonic---over the semidefinite cone---in the measure, and any
PSD matrix $N \in \symm_+^m$ can be decomposed as at most $m$ rank-one matrices.

\paragraph{Recovering the original inequalities.} Under the sufficient
conditions of~\S\ref{sec:sufficient-condition}, we can recover the original dual
problem by substituting $N = P_i^T N_i P_i$ for $i=0, \dots, d$.  This gives us
$d$ separate families of inequalities (one for each $i$) over both the smaller
matrices $N_i \in \symm_+^{m_i}$ and the nonnegative weights $w$. Namely:
\[
(b - A_0z)^T P_i^T N_i P_i (b - A_0z) \le \frac{z^TA_i^T P_i^T N_i P_i A_iz}{w_i}, ~~ \text{for all $N_i \in \symm_+^{m_i}$ and $w \ge 0$ with $\ones^Tw = 1$}.
\]
for each $i=1, \dots, d$, since $P_iA_j = 0$ for $i \ne j$. (Since the right
hand side is $0$ for $i=0$, we recover the original affine constraint as well.)
Finally, noting that this is true for $w = e_i$, the $i$th standard basis
vector, for $i=1, \dots, d$ gives the original
inequalities~\eqref{eq:psd-inequality}.

\paragraph{A (rewritten) problem.} Similar to problem~\eqref{eq:final-problem},
we can write
\begin{equation}\label{eq:final-problem-3}
\begin{aligned}
    & \text{minimize} && f(z)\\
    & \text{subject to} && (b - A_0z)^T N (b - A_0z) \le \sum_{i=1}^d \frac{z^TA_i^T N A_iz}{w_i},\\
    &&& ~~ \text{for all $N \in \symm_+^m$ and $w \ge 0$ with $\ones^Tw = 1$}.
\end{aligned}
\end{equation}
We can, of course, apply the same technique as before when $f$ is a quadratic to
get a dual problem that is a semidefinite program, with a family of constraints
over the nonnegative $w$. I suspect, though I do not prove, that picking a small
set of $w$ for this particular problem would yield reasonably tight lower bounds
on the optimal value of the original. On the other hand, it is worth noting that
most physics problems generally satisfy the sufficient conditions
of~\S\ref{sec:sufficient-condition}, so the original dual problem is likely to
be more useful in practice, since it also yields a much smaller total number of
variables.

\bibliographystyle{alpha}
\bibliography{citations.bib}

\end{document}

%% file: defs.tex
\newcommand{\ones}{\mathbf 1}
\newcommand{\reals}{{\mbox{\bf R}}}

\newcommand{\symm}{{\mbox{\bf S}}}  

\newcommand{\nullspace}{{\mathcal N}}

\newcommand{\cf}{{\it cf.}}
\newcommand{\eg}{{\it e.g.}}
\newcommand{\ie}{{\it i.e.}}

\newcommand{\BEAS}{\begin{eqnarray*}}
\newcommand{\EEAS}{\end{eqnarray*}}
\newcommand{\BEA}{\begin{eqnarray}}
\newcommand{\EEA}{\end{eqnarray}}
\newcommand{\BEQ}{\begin{equation}}
\newcommand{\EEQ}{\end{equation}}
\newcommand{\BIT}{\begin{itemize}}
\newcommand{\EIT}{\end{itemize}}

%% file: power-bound-generalization.bbl
\newcommand{\etalchar}[1]{$^{#1}$}
\begin{thebibliography}{MPHR{\etalchar{+}}16}

\bibitem[ADVB22]{angerisBoundsEfficiencyMetrics2022}
Guillermo Angeris, Theo Diamandis, Jelena Vu{\v c}kovi{\'c}, and Stephen Boyd.
\newblock Bounds on {{Efficiency Metrics}} in {{Photonics}}, April 2022.

\bibitem[AVB19]{angerisComputationalBoundsPhotonic2019}
Guillermo Angeris, Jelena Vu{\v c}kovi{\'c}, and Stephen Boyd.
\newblock Computational bounds for photonic design.
\newblock {\em ACS Photonics}, 6(5):1232--1239, May 2019.

\bibitem[AVB21a]{angerisConvexRestrictionsPhysical2021}
Guillermo Angeris, Jelena Vu{\v c}kovi{\'c}, and Stephen Boyd.
\newblock Convex restrictions in physical design.
\newblock {\em Scientific Reports}, 11(1):12976, December 2021.

\bibitem[AVB21b]{angerisHeuristicMethodsPerformance2021}
Guillermo Angeris, Jelena Vu{\v c}kovi{\'c}, and Stephen Boyd.
\newblock Heuristic methods and performance bounds for photonic design.
\newblock {\em Optics Express}, 29(2):2827, January 2021.

\bibitem[BV04]{cvxbook}
Stephen Boyd and Lieven Vandenberghe.
\newblock {\em Convex Optimization}.
\newblock {Cambridge University Press}, {Cambridge, United Kingdom}, first edition, 2004.

\bibitem[CSKD{\etalchar{+}}22]{chaoPhysicalLimitsElectromagnetism2022}
Pengning Chao, Benjamin Strekha, Rodrick Kuate~Defo, Sean Molesky, and Alejandro~W. Rodriguez.
\newblock Physical limits in electromagnetism.
\newblock {\em Nature Reviews Physics}, 4(8):543--559, August 2022.

\bibitem[DVY{\etalchar{+}}19]{doryInversedesignedDiamondPhotonics2019}
Constantin Dory, Dries Vercruysse, Kiyoul Yang, Neil Sapra, Alison Rugar, Shuo Sun, Daniil Lukin, Alexander Piggott, Jingyuan Zhang, Marina Radulaski, Konstantinos Lagoudakis, Logan Su, and Jelena Vu{\v c}kovi{\'c}.
\newblock Inverse-designed diamond photonics.
\newblock {\em Nature Communications}, 10(1):3309, December 2019.

\bibitem[GSJC20]{gustafssonUpperBoundsAbsorption2020}
Mats Gustafsson, Kurt Schab, Lukas Jelinek, and Miloslav Capek.
\newblock Upper bounds on absorption and scattering.
\newblock {\em New Journal of Physics}, March 2020.

\bibitem[Kes20]{kesiciMethodAntennaImpedance2020}
Cem Kesici.
\newblock Towards a method for antenna impedance tailoring through shape-optimization, with a bound on the associated cost-function.
\newblock Master's thesis, KTH Royal Institute of Technology, {Stockholm, Sweden}, 2020.

\bibitem[KM20]{kuangComputationalBoundsLight2020}
Zeyu Kuang and Owen~D. Miller.
\newblock Computational {{Bounds}} to {{Light}}\textendash{{Matter Interactions}} via {{Local Conservation Laws}}.
\newblock {\em Physical Review Letters}, 125(26):263607, December 2020.

\bibitem[LV13]{luNanophotonicComputationalDesign2013}
Jesse Lu and Jelena Vu{\v c}kovi{\'c}.
\newblock Nanophotonic computational design.
\newblock {\em Optics Express}, 21(11):13351, June 2013.

\bibitem[MCJR20]{moleskyGlobalOperatorBounds2020}
Sean Molesky, Pengning Chao, Weiliang Jin, and Alejandro Rodriguez.
\newblock Global {{T}} operator bounds on electromagnetic scattering: {{Upper}} bounds on far-field cross sections.
\newblock {\em Physical Review Research}, 2(3):033172, July 2020.

\bibitem[MCR20]{moleskyTOperatorLimitsElectromagnetic2020}
Sean Molesky, Pengning Chao, and Alejandro Rodriguez.
\newblock T-{{Operator Limits}} on {{Electromagnetic Scattering}}: {{Bounds}} on {{Extinguished}}, {{Absorbed}}, and {{Scattered Power}} from {{Arbitrary Sources}}.
\newblock {\em arXiv:2001.11531 [physics]}, January 2020.

\bibitem[MCR21]{moleskySionMinimaxTheorem2021}
Sean Molesky, Pengning Chao, and Alejandro~W. Rodriguez.
\newblock On {{Sion}}'s {{Minimax Theorem}}, {{Compact QCQPs}}, and {{Wave Scattering Optimization}}, May 2021.

\bibitem[MLF{\etalchar{+}}14]{menRobustTopologyOptimization2014}
H.~Men, K.~Lee, R.~Freund, J.~Peraire, and S.~Johnson.
\newblock Robust topology optimization of three-dimensional photonic-crystal band-gap structures.
\newblock {\em Optics Express}, 22(19):22632, September 2014.

\bibitem[MLP{\etalchar{+}}18]{moleskyInverseDesignNanophotonics2018}
Sean Molesky, Zin Lin, Alexander Piggott, Weiliang Jin, Jelena Vu{\v c}kovi{\'c}, and Alejandro Rodriguez.
\newblock Inverse design in nanophotonics.
\newblock {\em Nature Photonics}, 12(11):659--670, November 2018.

\bibitem[MPHR{\etalchar{+}}16]{millerFundamentalLimitsOptical2016}
Owen Miller, Athanasios Polimeridis, M.~T. Homer~Reid, Chiawei Hsu, {Brendan DeLacy}, John~D. Joannopoulos, Marin Solja{\v c}i{\'c}, and Steven~G. Johnson.
\newblock Fundamental limits to optical response in absorptive systems.
\newblock {\em Optics Express}, 24(4):3329, February 2016.

\bibitem[PLL{\etalchar{+}}15]{piggottInverseDesignDemonstration2015}
Alexander Piggott, Jesse Lu, Konstantinos Lagoudakis, Jan Petykiewicz, Thomas Babinec, and Jelena Vu{\v c}kovi{\'c}.
\newblock Inverse design and demonstration of a compact and broadband on-chip wavelength demultiplexer.
\newblock {\em Nature Photonics}, 9(6):374--377, June 2015.

\bibitem[SKLM21]{shimFundamentalLimitsMultifunctional2021}
Hyungki Shim, Zeyu Kuang, Zin Lin, and Owen~D. Miller.
\newblock Fundamental limits to multi-functional and tunable nanophotonic response, December 2021.

\bibitem[TAS{\etalchar{+}}20]{trivediBoundsScatteringAbsorptionless2020}
Rahul Trivedi, Guillermo Angeris, Logan Su, Stephen Boyd, Shanhui Fan, and Jelena Vu{\v c}kovi{\'c}.
\newblock Bounds for {{Scattering}} from {{Absorptionless Electromagnetic Structures}}.
\newblock {\em Physical Review Applied}, 14(1):014025, July 2020.

\bibitem[ZZM20]{zhaoMinimumDielectricResonatorMode2020}
Qingqing Zhao, Lang Zhang, and Owen~D. Miller.
\newblock Minimum {{Dielectric-Resonator Mode Volumes}}, August 2020.

\end{thebibliography}
